\documentclass[11pt,a4paper]{article}
\usepackage[utf8]{inputenc}

\input{macros}

\begin{document}

\title{Lyapunov stabilization of a nonlocal LWR traffic flow model}
\author{Jan Friedrich\footnotemark[1]}
\footnotetext[1]{RWTH Aachen University, Institute of Applied Mathematics, 52064 Aachen, Germany (friedrich@igpm.rwth-aachen.de).}
\date{\today}

\maketitle

\begin{abstract}
Using a nonlocal macroscopic LWR-type traffic flow model, we present an approach to control the nonlocal velocity towards a given equilibrium velocity. Therefore, we present a Lyapunov function measuring the $L^2$ distance between these velocities. We compute the explicit rate at which the system tends towards the stationary speed. The traffic is controlled by a leading vehicle. Numerical examples demonstrate the theoretical results and possible extensions of them.
\end{abstract}

  \medskip

  \noindent\textit{Keywords:} 
Lyapunov stabilization, nonlocal models, LWR model, macroscopic traffic flow


  \medskip

  \noindent\textit{AMS Subject Classification:} 
35L45, 35L65, 93D05
\section{Introduction}
In order to deal with the challenges arising from the progress in autonomous driving, classical approaches such as the Lighthill-Whitham-Richards (LWR) model \cite{lighthill1955kinematic, richards1956shockwaves} have been extended to include more information on the surrounding traffic, see for example~\cite{BlandinGoatin2016,friedrich2018godunov,KeimerPflug2017}.
These models are {nonlocal} traffic flow models.
Here, the flux function depends on an integral evaluation of the density or velocity.
In case of autonomous vehicles the integration area allows for an interpretation as a connection radius.

Nonlocal traffic flow models have been studied in various research directions over the years, e.g. well-posedness~\cite{chiarello2018global,friedrich2018godunov,KeimerPflug2017}, numerical schemes~\cite{ ChalonsGoatinVillada2018,friedrich2019maximum,friedrich2018godunov}, its singular limit behavior \cite{coclite2020general,colombo2019singular,keimer2019local} or modeling extensions~\cite{chiarello2019multiclass, keimer2019nonlocal, bayen2022multilane,friedrich2020nonlocal,friedrich2020onetoone}.
Nevertheless, there are only a few works concerning control problems \cite{bayen2021boundary, friedrich2022lyapunov, huang2020stability, karafyllis2022control}.
The works \cite{huang2020stability, karafyllis2022control} consider Lyapunov functions on a ring road.
In \cite{bayen2021boundary} the authors prove the exact controllability towards a target state together with explicit rates of convergence.
Instead of a ring road they consider a bounded domain and apply the control at the entrance and exit point of the road.
Recently, \cite{friedrich2022lyapunov} established Lyapunov stabilization of a second order traffic model on the micro- and macroscopic scale by controlling the leading vehicle.

Here, we consider as in \cite{bayen2021boundary,huang2020stability} the nonlocal LWR model of \cite{BlandinGoatin2016} for general velocity functions and want to obtain a Lyapunov function and an exponential decay rate.
In particular, in contrast to \cite{bayen2021boundary,huang2020stability}, the Lyapunov function considers the $L^2$ distance in the nonlocal velocity and we apply, as in \cite{friedrich2022lyapunov}, the control on the leading vehicle.

The upcoming work is structured as follows: Section \ref{sec:models} introduces the traffic flow model of \cite{BlandinGoatin2016}. Section \ref{sec:stable} contains our main result with the explicit rate on the decay rate. In the last section we present numerical examples which demonstrate the theoretical results.
Further, the numerical results probably hold under less restrictive conditions than we assume in the theory and can also be obtained for the microscopic scale.

\section{A nonlocal LWR traffic flow model}\label{sec:models}
We base our analysis of the Lyapunov stabilization on the following first order nonlocal traffic flow model:
\begin{align}
\label{eq:ND}
&\partial_t \rho(t,x)+\partial_x \left(\rho V(t,x)\right)=0,\qquad (t,x)\in \R_+\times \R\\
&V(t,x):=v\left(\wt \ast \rho\right)=v\left(\int_x^{x+\ndt} \wt(y-x)\rho(t,y)dy\right).\label{eq:ND:conv}
\end{align}
Here, $\rho$ is the traffic density, $v$ a suitable velocity function, $\ndt>0$ the nonlocal reach and $\wt$ a suitable kernel function.
 
This model was introduced in \cite{BlandinGoatin2016} and the well-posedness results were further extended in \cite{chiarello2018global,KeimerPflug2017}.
In particular, \cite{KeimerPflug2017} shows that no entropy condition is necessary as weak solutions are already unique.

The model \eqref{eq:ND:conv} is accompanied by initial conditions $\rho_0\in \BV(\R;[0,\rho_{\max}])$ and we impose the following assumptions on the kernel function and the velocity function:
\begin{assumption}\label{ass:vprime}
Let $v\in C^2([0,\rho_{\max}];\R_+)$ with $v'<0$ and in addition we denote the upper bound on the derivative of $v$ by $v'_{\max}$, i.e.
        \[ v'(\rho)\leq v'_{\max}<0\qquad \forall \ \rho\in[0,\rho_{\max}].\]
\end{assumption}
\begin{assumption}\label{ass:kernelcons}
We set the kernel function to
\begin{equation}
\wt(x)=\frac{1}{\ndt}.
\end{equation}
\end{assumption}
A constant convolution kernel can model, e.g., connected autonomous vehicles which have the same degree of accuracy on information about the downstream traffic, independent of the distance.
\begin{remark}
For the theoretical result we need to restrict ourselves to constant kernel functions.
In general, the kernel needs to fulfill
\begin{align} \label{eq:ass:kernel}
\wt \in C^1([0,\ndt];\R^+)\  & \text{with} \ \wt'\leq 0\ \text{and}\ \int_0^\ndt \wt(x) \dx=1 \ \forall \ndt>0
\end{align}
to guarantee existence and uniqueness, see \cite{BlandinGoatin2016, chiarello2018global}.
\end{remark}
In order to control the velocity to a stationary solution, i.e. $V(t,x)=\bar v$, for $t$ sufficiently large, we apply specific initial data:
\begin{assumption}\label{ass:initcond}
We consider strictly positive initial data $\rho_0\in \BV(\R;(0,\rho_{\max}])$ with
\begin{align}\label{eq:initcondcontrol}
\rho_0(x)=\bar \rho \qquad\text{for }x \geq b,
\end{align}
 $b\in\R$ and  $\bar \rho=v^{-1}(\bar v)$ for a given equilibrium velocity $\bar v$.
\end{assumption}
Note that $\bar \rho$ needs to be strictly positive such that $\bar v\in [0,v(0))$ holds.

As discussed in \cite[Remark 4.5]{KeimerPflug2017} the well-posedness of \eqref{eq:ND} with initial data $L^\infty(\R)\ni \rho_0 \not\in L^1(\R)$ is guaranteed.
In particular, the solution satisfies
\begin{align}\label{eq:maxprinciple}
0<\rho_{\min}:=\inf_{x\in\R}\rho_0(x)\leq \rho(t,x)\leq \sup_{x\in\R}\rho_0(x)\leq \rho_{\max},
\end{align}
see \cite{BlandinGoatin2016,chiarello2018global}, and hence we can set $v_{\max}'=v'(\rho_{\min})$.

\section{Lyapunov stabilization}\label{sec:stable}
The aim of this work is to control the nonlocal velocity $V(t,x)$ on the road towards a steady state $\bar v$.

In the following, we view the problem \eqref{eq:ND} with initial conditions fulfilling assumption \ref{ass:initcond} as an initial boundary value problem on the semi-infinite interval $(-\infty,\beta(t)]$ with $\beta(t):=b+\bar v t$ and we set the density for $x\geq\beta(t)$ to $\bar \rho$.
This induces the same waves as the initial value problem fulfilling assumption \ref{ass:initcond}, since for $x\geq\beta(t)$ the solution stays constant and $\beta(t)$ is the characteristic curve originating in $(0,b)$ in the $x$-$t$-plane, see \cite{KeimerPflug2017} for further details on the characteristics. 
This allows to view the left boundary as a leading vehicle which controls the traffic and moves with speed $\bar v$, see \cite{friedrich2022lyapunov} for further details.
This leading vehicle has a direct influence on the interval of length $\ndt$ behind it, i.e. $[\beta(t)-\ndt,\beta(t)]$.
Hence, our main result concerns this interval:
\begin{proposition}\label{prop:macroLyapunovDensity}
Let $0\leq \bar v<v(0)$ and the assumptions \ref{ass:vprime}, \ref{ass:kernelcons} and \ref{ass:initcond} hold.
We define
\begin{equation}\label{eq:macroLyapunov}
    L(t):=\int_{\beta(t)-\ndt}^{\beta(t)}(V(t,x)-\bar v)^2 dx.
\end{equation}
Then, we obtain the following bound
\begin{align*}
    L(t)\leq L(0) \exp\left(\frac{2}{\ndt} v'_{\max} \rho_{\min}t\right).
\end{align*}
\end{proposition}
\begin{proof}
First note that the boundary condition for $\beta(t)$ simplifies the nonlocal term in $\eqref{eq:ND:conv}$ and its partial derivatives for $x\in [\beta(t)-\ndt,\beta(t)]$ and a constant kernel as in the assumption \ref{ass:kernelcons}.
For simplicity, we define:
\begin{align*}
R(t,x):=\frac{1}{\ndt}\left( \int_x^{\beta(t)} \rho(t,y)\dy+(x+\ndt-\beta(t))\bar \rho\right).
\end{align*}
Next, we have
\begin{align}
    V(t,x)&=v\left(R(t,x)\right)\notag\\
    \partial_x V(t,x)&= \frac{1}{\ndt}\left(\bar \rho- \rho(t,x) \right)v'\left(R(t,x)\right)\label{eq:partialxV}
\end{align}
and hence, a direct relationship between the local density and the equilibrium density is obtained.
This derivative will play a key role for proving our stabilization result.
Next, we exploit the temporal derivative of the nonlocal term:
\begin{align}\notag
 \partial_t V(t,x)&=\frac1\ndt\left(\bar v \bar \rho+\int_x^{\beta(t)}\partial_t \rho(t,y)\dy-\bar v \bar \rho\right)v'\left(R(t,x)\right)\\
\notag &=-\frac{1}{\ndt}\int_x^{\beta(t)}\partial_x \left(\rho(t,y)V(t,y)\right)\dy v'\left(R(t,x)\right)\\
\label{eq:partialtV} &=\frac{1}{\ndt}\left(\rho(t,x)V(t,x)-\bar v\bar\rho \right) v'\left(R(t,x)\right).
\end{align}
Here, we used $\partial_t \rho(t,x)=-\partial_x (\rho(t,x)V(t,x))$.
Since $\rho$ is bounded and the integral evaluations are on a finite interval, we obtain that $R,\ V,\ \partial_x V,\ \partial_t V$ are all bounded.
In particular, this ensures that the integrals we consider in the following and the Lyapunov function $L(t)$ are bounded, too.

Turning now to the Lyapunov function \eqref{eq:macroLyapunov}, we directly obtain 
\begin{align*}
    \frac{d}{dt}L(t)=&\bar v (V(t,\beta(t))-\bar v)^2 - \bar v(V(t,\beta(t)-\ndt)-\bar v)^2\\
    &+2\int_{\beta(t)-\ndt}^{\beta(t)}(V(t,x)-\bar v)\partial_t V(t,x) dx.
\end{align*}
Due to the boundary conditions, we have $V(t,\beta(t))=\bar v$ and the first term equals zero.
Next, by using \eqref{eq:partialtV} and adding a zero we obtain
\begin{align*}
    \frac{d}{dt}L(t)=& - \bar v(V(t,\beta(t)-\ndt)-\bar v)^2+2\int_{\beta(t)-\ndt}^{\beta(t)}\frac{ v'\left(R(t,x)\right)}{\ndt}\rho(t,x)(V(t,x)-\bar v)^2 dx\\
    &+2\int_{\beta(t)-\ndt}^{\beta(t)}\frac{ v'\left(R(t,x)\right)}{\ndt}\bar v \left(\rho(t,x)-\bar\rho \right)(V(t,x)-\bar v)  dx.
\end{align*}
Using \eqref{eq:partialxV} we can rewrite the last term as
\begin{align*}
 &2\int_{\beta(t)-\ndt}^{\beta(t)}\frac{ v'\left(R(t,x)\right)}{\ndt}\bar v \left(\rho(t,x)-\bar\rho \right)(V(t,x)-\bar v)  dx \\
& =-\bar v \int_{\beta(t)-\ndt}^{\beta(t)}\partial_x (V(t,x)-\bar v)^2  dx   =\bar v\left(V(t,\beta(t)-\ndt)-\bar v\right)^2.
\end{align*}
Hence, we obtain
\begin{align*}
        \frac{d}{dt}L(t)=2\int_{\beta(t)-\ndt}^{\beta(t)}\frac{ v'\left(R(t,x)\right)}{\ndt}\rho(t,x)(V(t,x)-\bar v)^2 dx
        \leq \frac{2}{\ndt}\rho_{\min}v_{\max}'L(t).
\end{align*}
The latter inequality holds due to assumption \ref{ass:vprime} and the maximum principle for $\rho$ given by \eqref{eq:maxprinciple}.

Applying Gr\"onwall's lemma yields the result.
\end{proof}

\begin{remark}
In the case of a linear velocity function a stabilization result concerning the $L^2$ distance of the density and the stationary density is given by \cite[Proposition 3.8]{friedrich2022lyapunov}.
In this case, the Lyapunov function decreases only on a specific subinterval of $[\beta(t)-\ndt,\beta(t)]$.
\end{remark}
\begin{remark}
In order to estimate $\rho(t,x)$ from below and obtain a decreasing Lyapunov function, we need to ensure $\rho_{\min}>0$.
Hence, we cannot prove the exponential stability for $\bar v=v(0)$.
\end{remark}
As already mentioned in the introduction, \cite{huang2020stability} considers the stability of \eqref{eq:ND:conv} on a ring road. 
In addition, they restrict themselves to a linear velocity function.
One of the main results in \cite{huang2020stability} is that for a constant kernel as in assumption \ref{ass:kernelcons}, specific initial data and a specific nonlocal reach $\ndt$ traveling waves can be created, such that no control is possible.
Nevertheless, this is no contradiction to our results as we do not consider a ring road and such phenomena do not occur.

\section{Numerical simulations}
\subsection{Numerical scheme}
In order to illustrate the theoretical results we will use the numerical scheme presented in \cite{chiarello2019multiclass,friedrich2018godunov}.
For the spatial step size $\Dx>0$ the numerical flux is given by
\[F_{j+\frac12}^n=V_j^n \rho_j^n,\]
where
\[V_j^n = v\left(\sum_{k=0}^{\lfloor\ndt/\Delta x\rfloor-1}\gamma_k\rho^n_{j+k+1}\right)\quad\text{ and }\quad\gamma_k=\int_{k\Delta x}^{(k+1)\Delta x} \wt(y-x)dy.\]
The complete scheme is then given by
\begin{equation}\label{eq:macroscheme}
     \rho_j^{n+1}= \rho_j^{n}-\frac{\Dt^n}{\Dx}\left(F^n_{j+\frac12}-F^n_{j-\frac12}\right).
\end{equation}
At the right boundary we apply ghost cells equal to $\bar \rho$.
As CFL condition we choose an adaptive step size control determined by the maximum nonlocal velocity, i.e. $\Delta t^n\leq \Delta x/\max_j V_j^n $. 
\subsection{Numerical example}
As a numerical example we consider the situation of a bulk of cars standing bumper-to-bumper and starting to move.
For example such situations occur when a red light turns green or a completely stopped traffic jam begins to resolve.
This leads to the initial condition $\rho_0(x)=1$ for $x\leq 0=b$.
The speed law is given by $v(\rho)=1-\rho$.
Without a control the leading vehicle would move at maximum speed $v(0)=1$.
Instead we want to stabilize the traffic towards the velocity $\bar v=0.5$ (and $\bar \rho=0.5$).
In the numerical simulations the space step size is given by $\Dx=5\cdot 10^{-3}$.
The constant $\ndt$ is $\ndt=1$.
We compare the numerical stabilization for the following kernels
\begin{align*}
\wtcon(x)=\frac{1}{\ndt}, \qquad \wtlin(x)=\frac{2(\ndt-x)}{\ndt^2}, \qquad \wtconcave(x)=\frac{3(\ndt^2-x^2)}{2\ndt^3}.
\end{align*}
The theoretical results only ensure the convergence of the Lyapunov function \eqref{eq:macroLyapunov} for the kernel $\wtcon$.
Nevertheless, the kernels $\wtlin$ and $\wtconcave$ fulfill \eqref{eq:ass:kernel} such that the well-posedness of the model \eqref{eq:ND} is satisfied.
Figure \ref{fig:lyapunov} displays the logarithm of the Lyapunov function \eqref{eq:macroLyapunov} for all kernels together with the theoretical upper bound of Proposition \ref{prop:macroLyapunovDensity}.
As can be seen, even though not covered by our theoretical results, the bound is valid for all considered kernels.

\begin{figure}
\setlength{\fwidth}{0.5\textwidth}
\centering
%
%
\definecolor{mycolor1}{rgb}{0.92900,0.69400,0.12500}%
\definecolor{mycolor2}{rgb}{0.49400,0.18400,0.55600}%
\begin{tikzpicture}

\begin{axis}[%
width=1.4\fwidth,
height=0.5\fwidth,
at={(0\fwidth,0\fwidth)},
scale only axis,
xmin=0,
xmax=20,
xlabel style={font=\color{white!15!black}},
xlabel={$t$},
ymin=-35,
ymax=0,
axis background/.style={fill=white},
axis x line*=bottom,
axis y line*=left,
legend style={legend cell align=left, align=left, draw=white!15!black, legend pos=south west}
]
\addplot [color=black, line width=1.5pt]
  table[row sep=crcr]{%
0	-2.50749779388584\\
20	-22.5074977938858\\
};
\addlegendentry{exp. bound}

\addplot [color=blue, dotted, line width=1.5pt]
  table[row sep=crcr]{%
0	-2.12398992873574\\
0.23	-2.65207015585194\\
0.43	-3.08076706551632\\
0.629999999999999	-3.48298137061107\\
0.850000000000001	-3.89931996191734\\
1.09	-4.32787371718067\\
1.35	-4.76789312921934\\
1.63	-5.21940215510803\\
1.95	-5.71319817053829\\
2.31	-6.24689192283063\\
2.72	-6.83348233863492\\
3.19	-7.48547673533977\\
3.75	-8.24218896462096\\
4.44	-9.1546485646795\\
5.37	-10.3642755693211\\
7.02	-12.4881450330601\\
9.48	-15.6584525310043\\
11.27	-17.9849755005724\\
13.08	-20.3564204288524\\
15.01	-22.9041945466136\\
17.11	-25.6956798206577\\
19.41	-28.7723945447714\\
20	-29.5644822323621\\
};
\addlegendentry{$W_\eta^{\text{conc.}}$}

\addplot [color=gray, dashdotted, line width=1.5pt]
  table[row sep=crcr]{%
0	-2.02906485896463\\
0.340000000000003	-2.81384455813355\\
0.560000000000002	-3.28903264618946\\
0.780000000000001	-3.73738185041775\\
1.02	-4.19975537337172\\
1.28	-4.67465256737484\\
1.56	-5.16192068689521\\
1.88	-5.69495526549544\\
2.24	-6.27164050575811\\
2.66	-6.92207289409777\\
3.16	-7.67475622334278\\
3.81	-8.63136158352249\\
4.8	-10.0654683604577\\
7.1	-13.3928737039592\\
8.29	-15.137947548563\\
9.46	-16.8741536282152\\
10.66	-18.6754900094441\\
11.92	-20.5876150453382\\
13.27	-22.6574130740556\\
14.71	-24.886393279725\\
16.26	-27.30700538283\\
17.93	-29.9365903091389\\
19.72	-32.7767137538848\\
20	-33.2228060019614\\
};
\addlegendentry{$W_\eta^{\text{lin.}}$}

\addplot [color=red, dashed, line width=1.5pt]
  table[row sep=crcr]{%
0	-2.50749779388584\\
0.18	-2.85368313019897\\
0.359999999999999	-3.17444469237619\\
0.559999999999999	-3.50734561105066\\
0.780000000000001	-3.85092222617773\\
1.02	-4.20435015982042\\
1.29	-4.58089595568127\\
1.59	-4.97871032988068\\
1.92	-5.39674939875182\\
2.29	-5.84654091963512\\
2.71	-6.33849261231494\\
3.19	-6.8823917194172\\
3.74	-7.48781428250674\\
4.39	-8.18569970447271\\
5.17	-9.00577612657325\\
6.14	-10.0084166398938\\
7.42	-11.3141979340055\\
9.28	-13.1939313808556\\
12.53	-16.4598080515007\\
20	-23.9479408653679\\
};
\addlegendentry{$W_\eta^{\text{const.}}$}

\end{axis}

\end{tikzpicture}%
\captionsetup{margin=4cc}
\caption{Logarithm of the Lyapunov function \eqref{eq:macroLyapunov} over time together with the exponential upper bound for different kernel functions.}
\label{fig:lyapunov}
\end{figure}
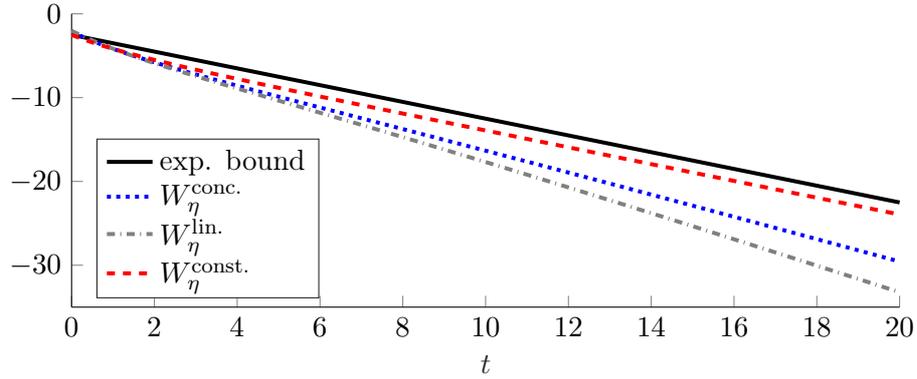
\begin{remark}
In order to prove the Lyapunov stabilization on the complete interval $[\beta(t)-\ndt,\beta(t)]$ we cannot consider the density, i.e.
\begin{equation}
\label{eq:lyapunovdensity}
\tilde{L}(t)=\int_{\beta(t)-\ndt}^{\beta(t)}(\rho(t,x)-\bar\rho)^2 dx.
\end{equation}
A numerical example with initial data
\[\rho_0(x)=\begin{cases}
0.01,\quad&\text{for }x\leq-0.5,\\
0.35,\quad&\text{for }x\in(-0.5,0],\\
0.5,\quad&\text{for }x>0,
\end{cases}\]
and the rest of the parameters as before can be seen in Figure \ref{fig:lyapunovdensity}.
In particular, the numerical Lyapunov function \eqref{eq:lyapunovdensity} has an increasing part which demonstrates that a stabilization result for \eqref{eq:lyapunovdensity} cannot be expected.
In order to get a decreasing Lyapunov function one needs to consider a specific subinterval of $[\beta(t)-\ndt,\beta(t)]$ as in \cite{friedrich2022lyapunov}.
\begin{figure}
\setlength{\fwidth}{0.5\textwidth}
\centering
%
%
\begin{tikzpicture}

\begin{axis}[%
width=1.4\fwidth,
height=0.5\fwidth,
at={(0\fwidth,0\fwidth)},
scale only axis,
xmin=0,
xmax=5,
xlabel style={font=\color{white!15!black}},
xlabel={$t$},
ymin=-2.1,
ymax=-2.02,
axis background/.style={fill=white},
axis x line*=bottom,
axis y line*=left,
legend style={legend cell align=left, align=left, draw=white!15!black, legend pos=south west}
]
\addplot [color=red, line width=1.5pt]
  table[row sep=crcr]{%
0	-2.03027049767338\\
0.0100000000000016	-2.03283456769106\\
0.0199999999999996	-2.03392741412071\\
0.0300000000000011	-2.03463345519142\\
0.0500000000000007	-2.03555656687742\\
0.0799999999999983	-2.03636258463457\\
0.120000000000001	-2.03692223981773\\
0.170000000000002	-2.03720354730551\\
0.239999999999998	-2.03719940327902\\
0.34	-2.03680539544987\\
0.859999999999999	-2.03418891775089\\
1	-2.0340276636442\\
1.13	-2.0341971990327\\
1.27	-2.03473107048799\\
1.47	-2.03592824293614\\
1.72	-2.03787943183065\\
1.91	-2.03966506040428\\
2.17	-2.04252572496358\\
2.45	-2.04610441044967\\
2.75	-2.05044975671868\\
3.08	-2.05575333999281\\
3.32	-2.05989490007671\\
3.58	-2.06460132781702\\
3.87	-2.07005484552805\\
4.59	-2.08418018967711\\
7.61	-2.14550295549081\\
8.26	-2.15884763127272\\
8.71	-2.16813033509133\\
9.4	-2.18246957187694\\
9.81	-2.19104812895909\\
10.5	-2.20561525730683\\
10.89	-2.21391510428396\\
11.56	-2.22831457688464\\
11.93	-2.236334517417\\
12.54	-2.24968792140892\\
12.89	-2.25741321958123\\
13.44	-2.2696707088145\\
13.77	-2.27708237578373\\
14.28	-2.28864548475882\\
14.59	-2.29572551806615\\
15.06	-2.30655853840569\\
15.37	-2.31375402516624\\
15.82	-2.32429545545883\\
16.11	-2.33113420440362\\
16.52	-2.34088856377457\\
16.81	-2.34783150567893\\
17.2	-2.35725113132765\\
17.47	-2.36381065369958\\
17.84	-2.37287734571835\\
18.11	-2.37953114473744\\
18.46	-2.38823064197875\\
18.71	-2.39447673855875\\
19.04	-2.40279081526513\\
19.29	-2.4091205541837\\
19.6	-2.41703495045507\\
19.83	-2.42293270487905\\
20	-2.42732395583838\\
};
\addlegendentry{\eqref{eq:lyapunovdensity}}

\end{axis}

\end{tikzpicture}%
\captionsetup{margin=4cc}
\caption{Logarithm of the Lyapunov function \eqref{eq:lyapunovdensity} over time with the kernel $\wtcon$.}
\label{fig:lyapunovdensity}
\end{figure}
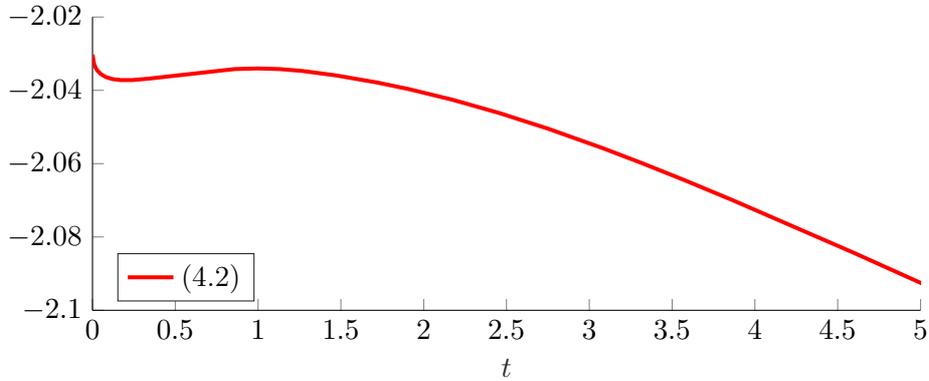
\end{remark}
\begin{remark}
In \cite{friedrich2022lyapunov} the Lyapunov stabilization of a corresponding microscopic model is derived with similar results as on the macroscopic scale.
For the model \eqref{eq:ND:conv} a microscopic formulation is considered in, e.g. \cite{ridder2019traveling}.
A numerical example of the situation above with $\rho_0(x)=1$ for $x\leq 0$ which considers the corresponding formulation of the Lyapunov function \eqref{eq:macroLyapunov} on the microscopic scale with $\wtcon$ can be seen in Figure \ref{fig:lyapunovmicro}.
The latter figure demonstrates that a stabilization to the equilibrium velocity can be obtained.
Nevertheless, for an analytic treatment one has to take care of discontinuities which can be seen in Figure \ref{fig:lyapunovmicro}.
These arise due to cars entering or leaving the domain $[\beta(t),\beta(t)-\ndt]$.  
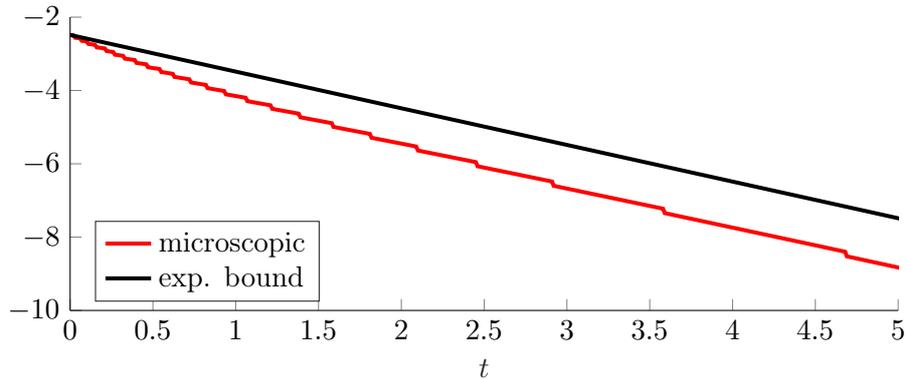
\begin{figure}
\setlength{\fwidth}{0.5\textwidth}
\centering
%
%
\begin{tikzpicture}

\begin{axis}[%
width=1.4\fwidth,
height=0.5\fwidth,
at={(0\fwidth,0\fwidth)},
scale only axis,
xmin=0,
xmax=5,
xlabel style={font=\color{white!15!black}},
xlabel={$t$},
ymin=-10,
ymax=-2,
axis background/.style={fill=white},
axis x line*=bottom,
axis y line*=left,
legend style={legend cell align=left, align=left, draw=white!15!black, legend pos=south west}
]
\addplot [color=red, line width=1.5pt]
  table[row sep=crcr]{%
0	-2.48500868539851\\
0.0199999999999996	-2.49511248731548\\
0.0300000000000011	-2.56146131714803\\
0.0599999999999987	-2.57711862788708\\
0.0700000000000003	-2.6449303061655\\
0.100000000000001	-2.66112715464796\\
0.109999999999999	-2.73045165788128\\
0.149999999999999	-2.75281527289702\\
0.16	-2.82373807197265\\
0.210000000000001	-2.85279474731768\\
0.219999999999999	-2.92539501924348\\
0.260000000000002	-2.94945479895265\\
0.27	-3.02375229426003\\
0.32	-3.05476997635049\\
0.329999999999998	-3.13086215269146\\
0.390000000000001	-3.16936266035147\\
0.399999999999999	-3.24730517805193\\
0.460000000000001	-3.28704818794978\\
0.469999999999999	-3.3669299252657\\
0.539999999999999	-3.41469835296032\\
0.550000000000001	-3.49663553212571\\
0.620000000000001	-3.54579248948708\\
0.629999999999999	-3.62984717779525\\
0.719999999999999	-3.69490045742919\\
0.73	-3.78117693818237\\
0.82	-3.84805912411096\\
0.830000000000002	-3.93668047465616\\
0.93	-4.01292901392363\\
0.940000000000001	-4.10400029511203\\
1.06	-4.19788615233364\\
1.07	-4.29154047592055\\
1.21	-4.40394543714121\\
1.22	-4.50032739275034\\
1.38	-4.63207018988976\\
1.39	-4.73133259679746\\
1.58	-4.89168533285892\\
1.59	-4.99398828012875\\
1.81	-5.18411860048867\\
1.82	-5.28964927605835\\
2.09	-5.52844333166061\\
2.1	-5.63740018214936\\
2.45	-5.95418159460614\\
2.46	-6.06679017124758\\
2.91	-6.4832879805378\\
2.92	-6.59982745039885\\
3.58	-7.22423542681243\\
3.59	-7.34494453064777\\
4.48	-8.20377284942187\\
4.68	-8.39892022848963\\
4.69	-8.52407010592367\\
5.62	-9.43668545444159\\
6.64	-10.4502102570702\\
7.79	-11.6063106128408\\
7.8	-11.7365538094065\\
8.17	-12.1051795807212\\
8.54	-12.4759386948177\\
9.02	-12.9747810804518\\
9.53	-13.5017750813916\\
10.17	-14.1378845293697\\
10.54	-14.5305600537015\\
11.56	-15.6250536510979\\
11.88	-15.9372341330726\\
12.22	-16.2852298947259\\
12.6	-16.6911335182207\\
13.51	-17.6722949483633\\
13.71	-17.8657452171832\\
14.04	-18.1953901160222\\
14.39	-18.5614336207656\\
14.84	-19.0502010867072\\
15.39	-19.6458158673632\\
15.59	-19.8466939262238\\
15.92	-20.1705387253628\\
16.26	-20.5205739754623\\
16.65	-20.9390940595585\\
17.47	-21.8254811559065\\
17.65	-22.001383892088\\
17.98	-22.3281821732148\\
18.33	-22.6915691154847\\
18.75	-23.1454871573651\\
19.41	-23.8607357936298\\
19.69	-24.1432711498926\\
20	-24.4543456319599\\
};
\addlegendentry{microscopic}

\addplot [color=black, line width=1.5pt]
  table[row sep=crcr]{%
0	-2.48500868539851\\
20	-22.4850086853985\\
};
\addlegendentry{exp. bound}

\end{axis}

\end{tikzpicture}%
\captionsetup{margin=4cc}
\caption{Logarithm of a microscopic equivalent to the Lyapunov function \eqref{eq:macroLyapunov} over time together with the corresponding exponential upper bound.}
\label{fig:lyapunovmicro}
\end{figure}
\end{remark}

\section{Conclusion}
In this work we have presented a Lyapunov function and the explicit rate such that the nonlocal velocity of a LWR-type traffic model on a single road tends to its equilibrium velocity.
For the theoretical analysis, we had to restrict ourselves to a constant kernel function.
Nevertheless, numerical examples suggest that the asymptotic stabilization effect can be obtained for general kernels and the corresponding microscopic model.
Future work may include extending the results to those cases.

Furthermore, by using the stabilization in the nonlocal velocity we are able to obtain our result on an interval of length $\ndt$ behind the leading vehicle. 
This is more intuitive than the stabilization result of \cite{friedrich2022lyapunov} which needs a specific subinterval.
Hence, future work may additionally include studying the Lyapunov stabilization of the nonlocal velocity for the second order traffic flow model considered in \cite{friedrich2020micromacro, friedrich2022lyapunov}.

\section*{Acknowledgement}
J.F. was supported by the German Research Foundation (DFG) under grant HE 5386/18-1, 19-2, 22-1, 23-1.

\vspace{\baselineskip}
\bibliographystyle{plain}
\providecommand{\othercit}{}
\providecommand{\jr}[1]{#1}
\providecommand{\etal}{~et~al.}

\end{document}